\documentclass[12pt]{article}
\usepackage{amsfonts}
\usepackage{amssymb}
\usepackage{mathrsfs}
\usepackage{srcltx}
\usepackage{dsfont,cases}
\usepackage{graphicx}
\usepackage{amsthm}
\usepackage{amstext}
\usepackage{amsopn}
\usepackage{graphicx}

\newtheorem{theorem}{Theorem}[section]
\newtheorem{lemma}[theorem]{Lemma}

\newtheorem{corollary}[theorem]{Corollary}

\theoremstyle{definition}
\newtheorem{definition}[theorem]{Definition}

\theoremstyle{remark}
\newtheorem{remark}[theorem]{Remark}

\newcommand{\ba}{\begin{array}}
\newcommand{\ea}{\end{array}}

\begin{document}
\date{}
\title{ \bf\large{Generalized Mountain Pass Lemma  Related with
a Closed Subset and Locally Lipschitz Functionals}}
\author{Fengying Li\textsuperscript{1}\footnote{Corresponding Author, Email: lify0308@163.com}\ \ Bingyu Li\textsuperscript{2}\ \ Shiqing Zhang\textsuperscript{3}
 \\
{\small \textsuperscript{1} School of Economic and Mathematics, Southwestern University of Finance and Economics\hfill{\ }}\\
\ \ {\small Chengdu, Sichuan, 611130, P.R.China.\hfill{\ }}\\
{\small \textsuperscript{2,3} Department of Mathematics, Sichuan University,\hfill{\ }}\\
\ \ {\small Chengdu, Sichuan, 610064, P.R.China.\hfill {\ }}}
\maketitle
\begin{abstract}
{The classical Mountain Pass Lemma of Ambrosetti-Rabinowitz  has been studied, extended and modified in several directions. Notable examples would certainly include the generalization to locally Lipschitz functionals by K.C. Chang, analyzing the structure of the critical set in the mountain pass theorem in the works of Hofer, Pucci-Serrin and Tian, and the extension by Ghoussoub-Preiss to closed subsets in a Banach space with recent variations. In this paper, we utilize the generalized gradient of Clarke and Ekeland's variatonal principle to generalize the Ghoussoub-Preiss's Theorem in the setting of locally Lipschitz functionals. We give an application to periodic solutions of Hamiltonian systems.}


\end{abstract}

\section{Introduction}

Saddle points in the Mountain pass Lemma (\cite{1}-\cite{23}) are different from maximum
points and minimum points. Maximum and Minimum problems in infinite dimensional space have a very long and prominent history (\cite{21}) with "isoperimetric problems" and the "problem of the brachistochrone" as two notable examples. In the 19th century ¡°Dirichlet principle¡± we essentially encountered the problem of minimizing a functional; however, complete rigor was mostly lacking and we had to wait for Hilbert for satisfactory completion of the Dirichlet principle. Continuing in the 20th century, Italian mathematician Tonelli introduced the concept of a weakly lower semi-continuous(w.l.s.c) functional and proved that a w.l.s.c functional defined on a weakly closed subset of a reflexive Banach space can attain its infimum if it is coercive(\cite{21}).
At times, the existence of a saddle point, which is neither a maximum nor minimum point, is of considerable importance. Minimax methods in the finite dimensional case(\cite{21}, \cite{23}) can be traced back to Birkhoff in 1917 and von Neumann's minimax theorem in 1928. We can also observe that the Mountain Pass Lemma of Ambrosseti-Rabinowitz(\cite{1}) in 1973 is a type of minimax theorem, which can be traced back to Courant in 1950 for the finite dimensional case(\cite{21}). From the finite dimensional case to the infinite dimensional case, the key step is using a Palais-Smale type compactness condition(PS) which implies Palais's Deformation Lemma. We should note the original proof of the Ambrosseti-Rabinowitz's Mountain Pass Lemma used Palais's Deformation Lemma(\cite{21}).
In the 1970's, Ekeland discovered a very important principle for lower semi-continuous functions on a complete metric space. Until the middle of the 1980's, Aubin-Ekeland(\cite{3}), Shi(\cite{20}) discovered the relationship between the Mountain Pass Lemma of Ambrosseti-Rabinowitz and Ekeland's variational priciple. The Mountain Pass Lemma of Ambrosseti-Rabinowitiz has been intensively studied and has found numerous applications (\cite{1}-\cite{23}). Of special note, it was generalized to the case of locally Lipschitz functionals by K.C. Chang(\cite{5}) where he also obtained more minimax theorems by using a deformation lemma.

In this paper, we use Ekeland's variational principle to prove a generalized Mountain Pass Lemma for locally Lipschitz functionals related with a closed subset, and we also found an applications to Hamiltonian systems with local Lispschtiz potential and a fixed energy.

\section{Classic Mountain Pass Theorem and Generalizations}
\setcounter{section}{1} \setcounter{equation}{0}

In 1973, Ambrosetti and Rabinowitz \cite{1} published the famous Mountain-Pass Theorem:
\begin{theorem}\label{thm 1}
{\rm (\cite{1})} Let $f$ be a $C^1-$real functional defined on a Banach space $X$ satisfying the following $(PS)$ condition:

Every sequence $\{x_n\}\subset X$ such that $\{f(x_n)\}$ is bounded and $\|f'(x_n)\|\rightarrow 0$ in $ X^*$ has a strongly convergent subsequence.

 Suppose there is an open neighborhood $\Omega$ of $x_0$ and a point $x_1\notin \bar{\Omega}$ such that
$
f(x_0),f(x_1)<c_0\leq\inf_{\partial\Omega}f,
$
and let
$$
\Gamma :=\{g\in C([0,1];X):g(0)=x_0, g(1)=x_1\}.
$$
Then
$
c :=\inf_{g\in\Gamma}\max_{t\in [0,1]}f(g(t))\geq c_0
$
 is a critical value of $f$: that is,
 there is $\bar{x}\in X$ such that $f(\bar{x})=c$ and $f'(\bar{x})=0$,
 where $f'(\bar{x})$ denotes the Fr\'{e}chet derivative of f at $\bar{x}$ .
\end{theorem}

Let $C^{1-0}(X;\mathds{R})$ be the space of locally Lipschitz mappings from
$X$ to $\mathds{R}$. For $\Phi\in C^{1-0}(X;\mathds{R})$ set (Clarke\cite{6})
$$\partial\Phi(x):=\{x^*\in X^*:<x^*,v>\leq\Phi^0(x,v),\forall v\in X\},$$
where
$\Phi^0(x,v):=\limsup\limits_{\mathop{{w\rightarrow x} }\limits_{t\downarrow 0}}\frac{\Phi(w+tv)-\Phi(w)}{t} $
denotes the generalized directional derivative of $\Phi$ at the point $x$ along the direction $v$. We should note that if $\Phi\in C^1$, then $\Phi^0(x,v)$ reduces to the G\^{a}teaux directional derivative and $\partial\Phi$ reduces to the classical derivative.

K.C. Chang \cite{5} generalized the classical (PS) condition and the Mountain Pass Theorem to local Lipschitz functions. Ribarska-Tsachev- Krastanov \cite{19} gave a generalization of a result of Chang for the case when "the separating mountain range has zero altitude" which is a version of the general mountain pass principle of Ghoussoub-Preiss for locally Lipschitz functions.

The generalization of the Mountain Pass Lemma of Ghoussoub-Preiss \cite{9} involves the modification of the classical Palais-Smale condition:
\begin{definition}\label{def1}
Let $X$ be a Banach space, $F$ a closed subset of $X$ and $\varphi$ a G\^{a}teaux-differentiable functional on $X$. The $(PS)_{F,c}$ condition is the following: if $\{x_n\}\subset X$ is a sequence satisfying the three conditions\\
(i) $d(x_n,F)\rightarrow 0.$\\
where $d(x,F):=\inf\limits_{y\in F}||x-y||$ denotes the distance between the point x and the set $F$,\\
(ii) $\varphi(x_n)\rightarrow c$,\\
(iii) $\varphi'(x_n)\rightarrow 0$,\\
then $\{x_n\}$ has a strongly convergent subsequence.
\end{definition}
\begin{definition}\label{def2}
Let $X$ be a Banach space, and $F\subset X$ a closed subset. We say $\Phi\in C^{1-0}(X;\mathds{R})$ meets the $(CPS)_{F,c}$ condition when the following is true: if $\{x_n\}\subset X$ satisfies
\begin{enumerate}
\item[(1)] $d(x_n,F)\rightarrow0$,
\item[(2)] $\Phi(x_n)\rightarrow c$,
\item[(3)] $(1+\|x_n\|)\cdot\min\limits_{y^*\in \partial\Phi(x_n)}\|y^*\|\rightarrow 0$,
\end{enumerate}
then $\{x_n\}$ has a convergent subsequence.
\end{definition}
We can define the $\delta$-distance (\cite{8}):
\begin{equation}\label{eq2}
\delta(x_1,x_2):=\inf\{l(c): c\in C^1([0,1],X),c(0)=x_1,c(1)=x_2\},
\end{equation}
where
\begin{equation}\label{eq3}
l(c):=\int^1_0\frac{\|\dot{c}(t)\|}{1+\|c(t)\|}dt.
\end{equation}
Then set $dist_{\delta}(x,F)=\inf\{\delta(x,y) : y\in F\}.$
\begin{definition}\label{def3}
Let $X$ be a Banach space, and $F\subset X$ a closed subset. We say $\Phi\in C^{1-0}(X;\mathds{R})$ meets the $(CPS)_{F,c;\delta}$ condition when the following is true: if $\{x_n\}\subset X$ satisfies
\begin{enumerate}
\item[(1)] $dist_{\delta}(x_n,F)\rightarrow0$,
\item[(2)] $\Phi(x_n)\rightarrow c$,
\item[(3)] $(1+\|x_n\|)\cdot\min\limits_{y^*\in \partial\Phi(x_n)}\|y^*\|\rightarrow 0$,
\end{enumerate}
then $\{x_n\}$ has a convergent subsequence.
\end{definition}
We can now state the Mountain-Pass Theorem generalized by Ghoussoub-Preiss\cite{9} for a continuous and G\^{a}teaux-differentiable functional statisfying the $(PS)_{F,c}$ condition:
\begin{theorem}{\rm(\cite{9})}\label{thm 2}
Let $\varphi: X\rightarrow \mathds{R}$ be a continuous and G\^{a}teaux-differentiable functional on a Banach space $X$ such that $\varphi':X\rightarrow X^*$ is continuous from the norm topology of $X$ to the $w^*-$topology of $X^*$.
Take $u,v\in X$, and let
$$
c:=\inf_{g\in\Gamma}\max_{t\in [0,1]}\varphi(g(t))
$$
 where $\Gamma=\Gamma_u^v$ is the set of all continuous paths joining $u$ and $v$. Suppose $F$ is a closed subset of $X$ such that $F\cap\{x\in X:\varphi(x)\geq c\}$ separates $u$ and $v$ and $\varphi$ satisfies the $(PS)_{F,c}$ condition, then there exists a critical point $\bar{x}\in F$ for $\varphi$ on $F$ with critical value $c$:
 $\varphi(\bar{x})=c,\varphi'(\bar{x})=0$.
\end{theorem}

In 2009, Goga \cite{10} studied a general Mountain Pass Theorem for local Lipschitz function. Let $(E,\|\cdot\|)$ be a Banach space, $S$ a compact metric space and $S_0$ a closed subset of $S$. Let $C(S,E)$ be the Banach space of all $E-$valued bounded continuous mapping on $S$ with the norm $\|\gamma\|:=\sup\limits_{x\in S} \|\gamma(x)\|$. Let $\gamma_0\in C(E,S)$ be a fixed element and define $$\Gamma=\{\gamma\in C(S,E):\gamma(s)=\gamma_0(s),\forall s\in S_0\},
c:=\inf_{\gamma\in \Gamma}\sup_{s\in S}f(\gamma(s)),$$ where $f$ is a real-valued function defined on $E$. Goga's result is the following:

\begin{theorem}\label{thm 3}
Let $f: E\rightarrow \mathds{R}$ be a locally Lipschitz function and $F$ a closed nonempty subset of $E$. Assume that
\begin{enumerate}
\item[(a)] $\gamma(S)\cap F\cap F_c\neq \emptyset, \forall \gamma\in \Gamma$, where $F_c=\{x\in E:f(x)\geq c\}$,
\item[(b)] $dist(\gamma_0(S_0),F)>0$, where $dist(\cdot,F)$ is the distance function to $F$ in $E$.
\end{enumerate}

Then for every $\varepsilon>0$ there exist $x_\varepsilon\in E$ such that
\begin{enumerate}
\item[(i)] $dist(x_\varepsilon,F)<\frac{3\varepsilon}{2}$,
\item[(ii)] $c\leq f(x_\varepsilon)<\varepsilon+\frac{5\varepsilon^2}{4}$,
\item[(iii)] $dist(0, \partial f(x_\varepsilon))\leq2\varepsilon$, where $\partial f(x)$ is the Clark sub-differential of $f$ at $x$.
\end{enumerate}

\end{theorem}

A key ingredient in the proof of Theorem \ref{thm 2} is provided by the following fundamental theorem in non-convex and nonlinear functional analysis established in the 1974 paper of Ivar Ekeland \cite{7}.

\begin{theorem}{\rm(\cite{7})}\label{thm 4} Let $(X,d)$ be a complete metric space with metric $d$ and $f: X\rightarrow \mathds{R}\cup\{+\infty\}$ a lower semi-continuous functional not identically $+\infty$ which is bounded from below. Let $\varepsilon>0$ and $u\in X$ such that
$
f(u)\leq\inf_{x\in X}f(x)+\varepsilon.
$
Then for any given $\lambda>0$, there exists $v_{\lambda}\in X$ such that $f(v_{\lambda})\leq f(u)$, $d(u,v_{\lambda})\leq\lambda$, and
$
f(w)>f(v_{\lambda})-\frac{\varepsilon}{\lambda} d(v_{\lambda},w),\ \ \forall w\neq v_{\lambda}.
$
\end{theorem}

Ekeland's variational principle has found numerous applications; in particular, we would like to observe that prior to Ghoussoub-Preiss\cite{9} it was used by Shi \cite{20} to prove a Mountain Pass Lemma and general min-max theorems for
locally Lipschitz functionals (K.C.Chang\cite{5}). In this paper, we will use Ekeland's variational principle to generalize the Ghoussoub-Preiss Theorem to the case of locally Lipschitz functional of
class $C^{1-0}$ satisfying the conditions $(CPS)_{F,c;\delta}$ or $(CPS)_{F,c}$.

\begin{theorem}\label{thm 5} Let $X$ be a Banach space with norm $||.||$, $C^{0}([0,1];X)$ the space of continuous mappings from $[0,1]$ to $X$, and $\Phi: X\rightarrow \mathds{R}$ a locally Lipschitz functional.
 For $z_0,z_1\in X$, define
\begin{equation}\label{eq4}
\begin{array}{l}
\Gamma :=\{c\in C^0([0,1];X):c(0)=z_0,c(1)=z_1\},
\gamma :=\inf\limits_{c\in \Gamma}\max\limits_{0\leq t\leq 1}\Phi(c(t)),
\end{array}
\end{equation}
and set
$
\Phi_{\gamma}:=\{x\in X:\Phi(x)\geq \gamma\}.
$
If $F\subset X$ is a closed subset such that $F\cap\Phi_{\gamma}$ separates $z_0$ and $z_1$, then there exists a sequence $\{x_n\}\subset X$ such that
$
dist_{\delta}(x_n,F)\rightarrow 0$, $\Phi(x_n)\rightarrow \gamma$ and
$(1+\|x_n\|)\min\limits_{y^*\in \partial\Phi(x_n)}\|y^*\|\rightarrow 0$.
\end{theorem}

\begin{theorem}\label{thm 6}
Under the assumptions of Theorem \ref{thm 5}, if we add that the set $F$ is
norm-bounded in the Banach space $X$, then a sequence
$\{x_n\}\subset X$ such that
$
d(x_n,F)\rightarrow 0$, $\Phi(x_n)\rightarrow \gamma$ and $(1+\|x_n\|)\min\limits_{y^*\in \partial\Phi(x_n)}\|y^*\|\rightarrow 0$.
\end{theorem}

\begin{theorem}\label{thm 7}
Under the assumptions of Theorem \ref{thm 5}, if $\Phi$ satisfies $(CPS)_{F,\gamma;\delta}$ condition, then $\gamma$ is a critical value for $\Phi$:
$\Phi(\bar{x})=\gamma,\ \ 0\in\Phi'(\bar{x})$.
\end{theorem}
\begin{theorem}\label{thm 8}
Under the assumptions of Theorem \ref{thm 5}, if we add the condition that the set $F$ is
bounded in the norm of the Banach space $X$, then we can change the
$(CPS)_{F,\gamma;\delta}$ condition to the $(CPS)_{F,\gamma}$ condition,
and conclude there exists a critical point $\bar{x}\in F$ for
$\Phi$ on $F$ with critical value \\
$\gamma: \Phi(\bar{x})=\gamma,0\in\Phi'(\bar{x})$.
\end{theorem}
\begin{remark}
The conclusions (i)-(iii) of Goga's Theorem \ref{thm 3} and the condition $(PS)_c$ in Ribarska- Tsachev- Krastanov \cite{19} are different from the conditions $(CPS)_{F,c}$ and $(CPS)_{F,c;\delta}$ stated here. Our results Theorem \ref{thm 5} and Theorem \ref{thm 6} are stronger since $(1+\|x_n\|)\min\limits_{y\ast\in \partial \Phi(x_n)}\|y\ast\|\rightarrow 0$ implies (iii) of the Theorem \ref{thm 3}. We would also like to note the assumptions in Theorem \ref{thm 7} and Theorem \ref{thm 8}, and our $(CPS)_{F,\gamma;\delta}$ and $(CPS)_{F,\gamma}$ conditions are weaker than those used \cite{10} and \cite{19}; therefore, the arguments in our paper differ from \cite{10} and \cite{19} since they could utilize the Borwein-Preiss variational principle or a deformation lemma, whereas we use the classical Ekeland's variational principle.
\end{remark}
\begin{remark}
We should note the difference between our Generalized Mountain Pass Lemma (GMPL) and the following theorem of Struwe(\cite{21}): Suppose $M$ is a closed convex subset of a Banach space $V$ and $E\in C^1(V)$ satisfies $(P.-S.)_M$ on $M$. Any sequence $\{u_n\}\subset M$ such that $|E(u_n)|\leq c$ uniformly, while\\ $g(u_m)=\sup\limits_{\mathop{{v\in M} }\limits_{\|u_m-v\|<1}}\langle u_m-v,DE(u_m)\rangle\rightarrow0$ $(m\rightarrow\infty)$, is relatively compact. Suppose further that $E$ admits two distinct relative minima $u_1$, $u_2$ in $M$. Then either $E(u_1)=E(u_2)=\beta$ and $u_1$, $u_2$ can be connected in any neighborhood of the set of relative minima $u\in M$ of $E$ with $E(u)=\beta$, or there exists a critical point $\bar{u}$ of $E$ in $M$ which is not a relative minimizer of $E$.

In Struwe's Theorem, $M$ is a closed convex subset of a Banach space, but in our GMPL we don't assume any convexity. We also don't assume that $E: M\rightarrow \mathds{R}$ possesses an extension $E\in C^1(V;\mathbb{R})$ to $V$, but only that the functional is locally Lipschitz. Struwe's Theorem assumes the existence of two local minimizers, but we only require the existence of two valleys which may not be local minimizers. In these ways, we see the premise in our GMPL is weaker than the corresponding conditions in Struwe's Theorem.
\end{remark}

\begin{remark}
The classical Mountain Pass Lemma and its many generalizations
are primarily concerned with ¡°saddle points¡±, but we should note the saddle
points encountered in these various Mountain Pass Lemmas are different from
those in the Von Neumann Minimax Theorem([23]). The Minimax Theorem
of Neumann is essentially related with convexity and concavity, whereas the
Mountain Pass Lemma is not related with convexity and concavity which is
related with (PS) compactness condiction and two valleys for functional. It
seems interesting to use Ekeland¡¯s variational principle to prove von Neumann Minimax Theorem.
\end{remark}

\section{The Proofs of Theorems \ref{thm 5}-\ref{thm 8}}
\setcounter{section}{1} \setcounter{equation}{0}

\begin{proof}:
Since the main ingredient is still Ekeland's variational principle, we utilize some notations and ideas from \cite{8} and \cite{9}, but must deviate in a few key steps.
Since the closed set $F_{\gamma}:=\Phi_{\gamma}\bigcap F$ separates $z_0$ and $z_1$, we can write
$
X\setminus F_{\gamma}:=\Omega_0\bigcup\Omega_1
$
where $z_0\in \Omega_0$, $z_1\in \Omega_1$ for open sets $\Omega_0$
and $\Omega_1$ with $\Omega_0 \cap \Omega_1=\emptyset$.

Choose $\varepsilon$ which satisfies
\begin{equation}\label{eq5}
0<\varepsilon<\frac{1}{2}\min\{1,dist_{\delta}(z_0,F_{\gamma}),dist_{\delta}(z_1,F_{\gamma})\}.
\end{equation}
By the definition of $\Gamma$, we can find $c\in\Gamma$ such that
\begin{equation}\label{eq6}
\max_{0\leq t\leq1}\Phi(c(t))<\gamma+\frac{\varepsilon^2}{4}.
\end{equation}
If we define $t_0$ and $t_1$ by
\begin{eqnarray*}
t_0&:=&\sup\{t\in[0,1]:c(t)\in\Omega_0,dist_{\delta}(c(t),F_{\gamma})\geq\varepsilon\},\\
t_1&:=&\inf\{t\in[t_0,1]:c(t)\in\Omega_1,dist_{\delta}(c(t),F_{\gamma})\geq\varepsilon\},
\end{eqnarray*}
then since $c(0)=z_0$, we have by (\ref{eq5}) and the continuity of c that
$t_0>0$; moreover, by $c(t_0)\in\bar{\Omega}_{0}$ and
$dist_{\delta}(c(t_0),F_{\gamma})\geq\epsilon$, we have
$c(t_0)\in\Omega_0$. Then $\Omega_0\cap\Omega_1=\emptyset$ implies
$t_1>t_0$. Again by (\ref{eq5}) and the continuity of $c$ we have $t_1<1$. So altogether
$0<t_0<t_1<1$.
Let
\begin{equation}\label{eq7}
\Gamma(t_0,t_1):=\{f\in C^0([t_0,t_1],X):f(t_0)=c(t_0),f(t_1)=c(t_1)\},
\end{equation}
and consider the following distance in $\Gamma(t_0,t_1)$:
\begin{equation}\label{eq8}
\rho(f_1,f_2):=\max_{t_0\leq t\leq t_1}\delta(f_1(t),f_2(t)),
\end{equation}
where
\begin{equation}\label{eq9}
\delta(x_1,x_2):=\inf\{l(c) : c\in C^1([0,1],X),c(0)=x_1,c(1)=x_2\}
\end{equation}
with
\begin{equation}\label{eq10}
l(c):=\int^1_0\frac{\|\dot{c}(t)\|}{1+\|c(t)\|}dt.
\end{equation}
For $x\in X$, we define the function
\begin{equation}\label{eq11}
\Psi(x):=\max\{0,\varepsilon^2-\varepsilon dist_{\delta}(x,F_{\gamma})\}.
\end{equation}
 A map $\varphi: \Gamma(t_0,t_1)\rightarrow \mathds{R}$ is defined by
 \begin{equation}\label{eq12}
\varphi(f):=\max_{t_0\leq t\leq t_1}\{\Phi(f(t))+\Psi(f(t))\}.
\end{equation}
Since $f(t_0)=c(t_0)\in\Omega_0, f(t_1)=c(t_1)\in\Omega_1$,
there exists $t_f\in (t_0,t_1)$ satisfying $f(t_f)\in \partial\Omega_0\subset F_{\gamma}$; therefore,
\begin{equation}\label{eq14}
dist_{\delta}(f(t_f),F_{\gamma})=0,
\end{equation}
and for $\forall f\in\Gamma(t_0,t_1)$, we have
\begin{equation}\label{eq15}
\varphi(f)\geq\Phi(f(t_f))+\Psi(f(t_f))\geq\gamma+\varepsilon^2.
\end{equation}
On the other hand, if we denote $\hat{c}=c|_{[t_0,t_1]}$, then
\begin{equation}\label{eq16}
\varphi(\hat{c})\leq \max_{0\leq t\leq 1}\{\Phi(c(t))+\Psi(c(t))\}\leq\gamma+\frac{5}{4}\varepsilon^2.
\end{equation}
Notice that $\Gamma(t_0,t_1)$ is a complete metric space \cite{7},\cite{8}.
Since $\Phi$ and $\Psi$ are lower semi-continuous, so is $\varphi$. Now (\ref{eq15}) implies $\phi$ has a lower bound, and by (\ref{eq15}) and (\ref{eq16}) we have
\begin{equation}\label{eq17}
\varphi(\hat{c})\leq\inf \varphi+\frac{\varepsilon^2}{4}.
\end{equation}
In Ekeland's variational principle, we use
$\frac{\varepsilon^2}{4}$ in place of $\epsilon$, and take
$\lambda=\frac{\epsilon}{2}$, then there exists $\hat{f}\in
\Gamma(t_0,t_1)$ such that
\begin{equation}\nonumber
\begin{array}{l}
\varphi(\hat{f})\leq\varphi(\hat{c}),\ \
\rho(\hat{f},\hat{c})\leq\frac{\varepsilon}{2},\ \
\varphi(f)\geq\varphi(\hat{f})-\frac{\varepsilon}{2}\rho(f,\hat{f}), \ \ \forall f\in \Gamma(t_0,t_1).\label{eq18}
\end{array}
\end{equation}
Let
\begin{equation}\label{eq19}
M:=\{t\in [t_0,t_1]:\Phi(\hat{f}(t))+\Psi(\hat{f}(t))=\varphi(\hat{f})\}.
\end{equation}
The claim is that $M$ is a non-empty compact set which avoids $t_0$ and $t_1$.

By the definitions of $t_0$ and $t_1$, we have
\begin{equation}\label{eq20}
dist_{\delta}(c(t_i),F_{\gamma})\geq \varepsilon, \ \ i=0,1,
\end{equation}
so $\Psi(\hat{c}(t_i))=0$. By (\ref{eq6}) and (\ref{eq15}) we have
\begin{equation}\label{eq21}
\Phi(\hat{f}(t_i))+\Psi(\hat{f}(t_i))\leq\Phi(\hat{c}(t_i))+\Psi(\hat{c}(t_i))\leq \gamma+\frac{\varepsilon^2}{4}<\varphi(\hat{f}),i=0,1,
\end{equation}
which implies $t_0,t_1\notin M$.
\vspace{0.2cm}

\textbf{\underline{Claim}}: There exists $t\in M$ such that
\begin{equation}\label{eq22}
\min_{x^*\in \partial\Phi(\hat{f}(t))}\|x^*\|(1+\|\hat{f}(t)\|)\leq \frac{3\varepsilon}{2};
\end{equation}

\begin{proof}: If not, for any $t\in M$,
\begin{equation}\label{eq23}
\min_{x^*\in \partial\Phi(\hat{f}(t))}\|x^*\|(1+\|\hat{f}(t)\|)>\frac{3\varepsilon}{2}.
\end{equation}
It is well known that $\|x^*\|=\sup_{v\neq 0}\frac{<x^*,v>}{\|v\|}$ where
\begin{equation}\label{eq24}
x^*\in \partial\Phi(\hat{f}(t))=\{x^*\in X^*:<x^*,v>\leq\Phi^0(\hat{f}(t),v),\forall v\in X\},
\end{equation}
where the following definition
$$\Phi^0(x,v)=\limsup_{\mathop{{w\rightarrow x} }\limits_{t\downarrow 0}}\frac{\Phi(w+tv)-\Phi(w)}{t} $$
denotes the generalized directional derivative of $\Phi$ at the point $x$ along the direction $v$.
Notice that
\begin{equation}\label{eq25}
\Phi^0(x,v)=\max\{<x^*,v>:x^*\in\partial\Phi(x)\}.
\end{equation}
Then for all $t\in M$ there exists $u(t)\in X$ such that $\|u(t)\|=(1+\|\hat{f}(t)\|)$ and
\begin{equation}\label{eq26}
\Phi^0(\hat{f}(t),u(t))<-\frac{3\varepsilon}{2}.
\end{equation}

Let $t\in M$ be such that (\ref{eq23}) holds, then
$\partial\Phi(\hat{f}(t))\bigcap\frac{3\epsilon}{2(1+||\hat{f}(t)||)}B_{X^{*}}=\emptyset.$
Notice that sets $\partial\Phi(\hat{f}(t))$ and
$\frac{3\epsilon}{2(1+||\hat{f}(t)||)}B_{X^{*}}$ are convex and
$w^*$ compact, so by the separation theorem, the two sets can be
separated by an element of $X$; that is, there is $v_0\in X$ such that
$||v_0||=1$ and
$$\sup\{<x^*,v_0>:x^*\in\frac{3\epsilon}{2(1+||\hat{f}(t)||)}B_{X^{*}}\}<\inf\{<x^*,v_0>:x^*\in\partial\Phi(\hat{f}(t))\}$$
Notice that the left side of the above inequality is just
$\frac{3\epsilon}{2(1+||\hat{f}(t)||)}$. Hence if we let
$h=-(1+||\hat{f}(t)||)v_0$, then we have
$$\sup\{<x^*,h>:x^*\in\partial\Phi(\hat{f}(t))\}<-\frac{3\epsilon}{2}.$$
Notice that the left side of the above inequality is equal to
$\Phi^0(\hat{f}(t),h)$, so we get (\ref{eq26}) for $u(t)=h$. Let
$N(t):=\{s\in M:\Phi^0(\hat{f}(s),u(t))<-\frac{3}{2}\epsilon\}.$
Since $x\mapsto\Phi^0(x,v)$ is upper semicontinuous for any given $v$, $N(t)$ is an open subset of $M$ and
$M$ can be covered by the open sets $N(t)$ for $t\in M.$
Since $M$ is compact, we can pick a finite open sub-cover for $M$, $\left\{ N(t_k):0\leq k\leq K \right\}$. Then for the partition of unity associated with this cover on $M$, there are continuous functions $\xi_k(t):0\leq\xi_k(t)\leq 1$ for $0\leq k\leq K$ with $\sum_1^K\xi_k(t)=1.$

Let $v(t):=\sum_{k=0}^K\xi_k(t)u(t_k)$ and observe the continuous map $v: M\rightarrow X$ satisfies
\begin{equation}\label{eq27}
\Phi^0(\hat{f}(t),v(t))<-\frac{3\varepsilon}{2}, \ \ \|v(t)\|\leq\sum_{k=0}^K\xi_k(t)(1+\|\hat{f}(t)\|=(1+\|\hat{f}(t)\|).
\end{equation}
Since $M\subset[t_0,t_1]$ and M is a nonempty compact set with
$t_0,t_1\notin M$, so by Tietze extension theorem, we can extend $v$
to a continuous function defined on $[t_0,t_1]$ (which we still
denote by $v$) which satisfies $v(t_0)=v(t_1)=0$ and
\begin{equation}\label{eq28}
 \|v(t)\|\leq(1+\|\hat{f}(t)\|), \forall t\in[t_0,t_1].
\end{equation}
Since $v(t_0)=v(t_1)=0$, $\forall h>0$, $\hat{f}+hv\in\Gamma$; hence,
\begin{equation}\label{eq29}
 \varphi(\hat{f}+hv)\geq \varphi(\hat{f})-\frac{\varepsilon}{2}\rho(\hat{f}+hv,\hat{f}).
\end{equation}
Choose $t_h\in[t_0,t_1]$ which satisfies:
\begin{equation}\label{eq30}
 \varphi(\hat{f}+hv)=(\Phi+\Psi)(\hat{f}(t_h)+hv(t_h)).
\end{equation}
Notice that here $t_h$ is defined for each $h>0$. By the definition
of $\varphi$, we know that for any $h>0$, there holds $\varphi(\hat{f})\geq (\Phi+\Psi)(\hat{f}(t_h))\label{eq31}$.

So $\forall h>0$ we have
\begin{equation}\label{eq32}
(\Phi+\Psi)(\hat{f}(t_h)+hv(t_h))\geq(\Phi+\Psi)(\hat{f}(t_h))-\frac{\varepsilon}{2}\rho(\hat{f}+hv,\hat{f});
\end{equation}
that is,
\begin{equation}\label{eq33}
\Phi(\hat{f}(t_h)+hv(t_h))-\Phi(\hat{f}(t_h))\geq-\Psi(\hat{f}(t_h)+hv(t_h))+\Psi(\hat{f}(t_h))-\frac{\varepsilon}{2}\rho(\hat{f}+hv,\hat{f})).
\end{equation}
If we recall the definition of $\Psi$, then $\Psi$ is $\varepsilon-$Lipschitz, and so the above inequality implies
\begin{equation}\label{eq34}
\Phi(\hat{f}(t_h)+hv(t_h))-\Phi(\hat{f}(t_h))\geq-\frac{3\varepsilon}{2}\rho(\hat{f}+hv,\hat{f})).
\end{equation}
Notice that if $h_n\rightarrow 0^+$, we can pass to a sequence $\{t_{h_n}\}$ with $t_{h_n}\rightarrow\tau\in M$ since $M$ is compact. Calculating
\begin{equation}\label{eq35}
\limsup_{n\rightarrow +\infty}\frac{\Phi(\hat{f}(t_{h_n})+h_nv(t_{h_n}))-\Phi(\hat{f}(t_{h_n}))}{h_n}
\geq-\frac{3\varepsilon}{2}\liminf_{n\rightarrow +\infty}\frac{\rho(\hat{f}+h_nv,\hat{f}))}{h_n}.
\end{equation}
and further by $\Phi\in C^{1-0}$ and the definitions of Clark's generalized gradient and the metric $\rho$,
 we have
\begin{equation}\label{eq36}
\Phi^0(\hat{f}(\tau),v(\tau))\geq-\frac{3\varepsilon}{2}\max_{t_0\leq t\leq t_1}(\frac{\|v(t)\|}{1+\|\hat{f}(t)\|})\geq-\frac{3\varepsilon}{2}.
\end{equation}
In fact, by $\Phi\in C^{1-0}$ and the continuity for $v(t)$, we know
that
$$\frac{\Phi(\hat{f}(t_{h_n})+h_nv(t_{h_n}))-\Phi(\hat{f}(t_{h_n})+h_nv(\tau))}{h_n}\leq L|v(t_{h_n})-v(\tau)|\rightarrow 0,$$
hence
\begin{eqnarray*}
&\lim&\sup\limits_{n\rightarrow +\infty}\frac{\Phi(\hat{f}(t_{h_n})+h_nv(t_{h_n}))-\Phi(\hat{f}(t_{h_n}))}{h_n}\\
&\leq& \limsup_{n\rightarrow +\infty}\frac{\Phi(\hat{f}(t_{h_n})+h_nv(t_{h_n}))-\Phi(\hat{f}(t_{h_n})+h_nv(\tau))}{h_n}\\
&+&\limsup_{n\rightarrow +\infty}\frac{\Phi(\hat{f}(t_{h_n})+h_nv(\tau))-\Phi(\hat{f}(t_{h_n})}{h_n}\\
&=&\Phi^0(\hat{f}(\tau),v(\tau))
\end{eqnarray*}
Using the definition (\ref{eq8}) of the metric $\rho$, we have
\begin{eqnarray*}
&\rho&(\hat{f}+h_nv,\hat{f}))\\
&=&\max_{t_0\leq t\leq t_1}\inf\{\int_0^1\frac{||\dot{c}(s)||}{1+||c(s)||}ds,c(s)\in C^0([0,1],X),c(0)=\hat{f},c(1)=\hat{f}+h_nv\}.
\end{eqnarray*}
Specifically, if we take the following loop connecting $\hat{f}$ and $\hat{f}+h_nv$:
$$c(s)=(1-s)\hat{f}(t)+s(\hat{f}(t)+h_nv(t)),$$
then we have $\dot{c}(s)=h_nv(t)$, and so
$$c(s)=(1-s)\hat{f}(t)+s(\hat{f}(t)+h_nv(t))\rightarrow \hat{f}(t),n\rightarrow +\infty.$$
So we have that
$$\liminf_{n\rightarrow +\infty}\frac{\rho(\hat{f}+h_nv,\hat{f}))}{h_n}\leq\max_{t_0\leq t\leq t_1}(\frac{\|v(t)\|}{1+\|\hat{f}(t)\|}),$$
and (\ref{eq36}) is proved, which violates (\ref{eq27}) and shows that we cannot
have the inequality (\ref{eq23}); therefore, there is $\bar{t}\in M$ such
that
\begin{equation}\label{eq37}
\min_{x^*\in \partial\Phi(\hat{f}(\bar{t}))}\|x^*\|(1+\|\hat{f}(\bar{t})\|)\leq\frac{3\varepsilon}{2}.
\end{equation}

\end{proof}
By the definitions of $t_0$ and $t_1$, we have that
$d(\hat{c}(t),F_{\gamma})\leq\epsilon$ for
$t_0<t<t_1$; furthermore, by continuity of $\hat{c}(t)$
and $d(x,F_{\gamma})$ on $x$, we have that
$$d(\hat{c}(t),F_{\gamma})\leq\epsilon,\forall t\in[t_0,t_1].$$
Notice that here $d(\hat{c}(t),F_{\gamma})$ is the distance between
$\hat{c}(t)$ and $F_{\gamma}$ deduced by the norm in the Banach
space $X$. We use the notation $dis_{\delta}(\hat{c}(t),F_{\gamma})$
to denote the distance between $\hat{c}(t)$ and $F_{\gamma}$ deduced
by $\delta$ in (\ref{eq9}). By the definitions of $\delta$ and the norm,
we have that
$$\delta(x_1,x_2)\leq ||x_1-x_2||$$
so
$$dis_{\delta}(\hat{c}(t),F_{\gamma})\leq d(\hat{c}(t),F_{\gamma})\leq\epsilon,\forall t\in[t_0,t_1].$$
We notice that $\rho$ is the distance deduced by $\delta$ in (\ref{eq9}). Since $\rho(\hat{f},\hat{c})\leq \frac{\varepsilon}{2}$, the triangle inequality implies that for all $t\in [t_0,t_1]$ we have
\begin{equation}\label{eq38}
dis_{\delta}(\hat{f}(t),F_{\gamma})\leq\frac{\varepsilon}{2}+dis_{\delta}(\hat{c}(t),F_{\gamma}))\leq \frac{\varepsilon}{2}+\varepsilon=\frac{3\varepsilon}{2}.
\end{equation}
Set $x=\hat{f}(\bar{t})$, we get $dis_{\delta}(x,F_{\gamma})\leq\frac{3\varepsilon}{2}.$

If $F$ is bounded and a closed subset of $X$, then by the definition
of $\delta$, we know that (\cite{7}) $\delta$ distance is equivalent to
the norm distance, so there is $c>0$ such that $dis_{\delta}(x,F_{\gamma})\geq c d(x,F_{\gamma}).$
Then $\varphi(\hat{f})\leq \varphi(\hat{c})$ yields
\begin{equation}\label{eq39}
\gamma+\varepsilon^2\leq\Phi(\hat{f}(\bar{t}))+\Psi(\hat{f}(\bar{t}))\leq\gamma+\frac{5\varepsilon^2}{4}.
\end{equation}
Then we get
\begin{equation}\nonumber
\begin{array}{l}
\min_{x^*\in \partial\Phi(x)}\|x^*\|(1+\|x\|)\leq\frac{3\varepsilon}{2},\label{eq40}\\
d(x,F_{\gamma})\leq\frac{1}{c}\frac{3\varepsilon}{2},\\
\gamma\leq\Phi(x)\leq\gamma+\frac{5\varepsilon^2}{4}.
\end{array}
\end{equation}
If we let $\varepsilon=\frac{1}{n}\rightarrow 0$, then we arrive at
a sequence $\{x_n\}$ which satisfies the requirements of Theorems
\ref{thm 5} and \ref{thm 6}. Theorems \ref{thm 7} and \ref{thm 8} follow from Theorems \ref{thm 5} and
\ref{thm 6}.

\end{proof}

\section{An Application to Hamiltonian systems}
\setcounter{section}{1} \setcounter{equation}{0}

Let $V\in C^{1-0}(\mathds{R}^n,\mathds{R})$; that is, $V$ is a locally Lipschitz potential function defined on $\mathds{R}^n$.
Let us consider the second order Hamiltonian systems
\begin{numcases}{}
-\ddot{q}(t)\in \partial V(q) \label{eq41}\\
\frac{1}{2}|\dot{q}|^2+V(q)=h\in \mathds{R}\label{eq42}
\end{numcases}
\begin{theorem}\label{thm 9}
Suppose $V\in C^{1-0}(\mathds{R}^n,\mathds{R})$ and $h\in \mathds{R}$
satisfy
  \begin{enumerate}
    \item[($V_1$)] $ V(-q)=V(q);$

    \item[($V_2$)] $\exists \mu_1>0,\mu_2\geq 0$, such that $\langle y,q\rangle\geq\mu_1V(q)-\mu_2,$ $\forall y\in \partial V(q)$, $\forall q\in \mathds{R}^n$;

    \item[($V_3$)] $V(q)\geq h, |q|\rightarrow +\infty.$
  \end{enumerate}
Then for any $h>\frac{\mu_2}{\mu_1},$ the system $(\ref{eq41})-(\ref{eq42})$ has at least
one non-constant periodic solution with the given energy $h$ which can be obtained by
Theorem \ref{thm 7}.
\end{theorem}

\begin{corollary}
For $a>0$, $\mu_1>0$, $\mu_2\geq0$, let $V(q)=a|q|^{\mu_1}+\frac{\mu_2}{\mu_1}$.
Then for any $h>\frac{\mu_2}{\mu_1},$ the system $(\ref{eq41})-(\ref{eq42})$ has at least
one non-constant periodic solution with the given energy $h$ which can be obtained by
Theorem \ref{thm 7}.
\end{corollary}
\begin{remark}If $\mu_1>1$, then $V(q)=a|q|^{\mu_1}+\frac{\mu_2}{\mu_1}\in C^1(\mathds{R}^n, \mathds{R})$; but if $0<\mu_1\leq1$ it is not in $C^1(\mathds{R}^n,\mathds{R})$, but $V\in C^{1-0}(\mathds{R}^n,\mathds{R})$.
\end{remark}
In order to prove Theorem \ref{thm 9}, we define the Sobolev space
\begin{equation}\label{eq43}
H^1:=W^{1,2}(\mathds{R}/T\mathds{Z},\mathds{R}^n)=\{u:\mathds{R}\rightarrow \mathds{R}^n,u\in L^2,\dot{u}\in L^2,u(t+1)=u(t)\}.
\end{equation}
 Then the standard $H^1$ norm is equivalent to
\begin{equation}\label{eq44}
 \|u\|:=\|u\|_{H^1}=\left(\int^1_0|\dot{u}|^2dt\right)^{1/2}+|\int_0^1 u(t)dt|.
 \end{equation}
\begin{lemma}\label{lemma1}{\rm(\cite{2})}\ \ Let
$f(u):=\frac{1}{2}\int^1_0|\dot{u}|^2dt\int^1_0(h-V(u))dt$ and
$\widetilde{u}\in H^1$ be such that $f^{\prime}(\widetilde{u})=0$
 and $f(\widetilde{u})>0$. Set
\begin{equation}\label{eq45}
\frac{1}{T^2}:=\frac{\int^1_0(h-V(\widetilde{u}))dt}{\frac{1}{2}\int^1_0|\dot{\widetilde{u}}|^2dt}.
\end{equation}
Then $\widetilde{q}(t)=\widetilde{u}(t/T)$ is a non-constant
$T$-periodic solution for {\rm(\ref{eq41})-(\ref{eq42})}.
\end{lemma}

In a manner similar to Ambrosetti-Coti
Zelati\cite{2}, from the symmetry condition $(V_1)$ we let
$E:=\{u\in
H^1=W^{1,2}(\mathds{R}/\mathds{Z},\mathds{R}^n),u(t+1/2)=-u(t)\}.$
A similar proof as in \cite{2}, we have
\begin{lemma}\label{lemma2}
 If $\bar{u}\in E$ is a critical point of
$f(u)$ and $f(\bar{u})>0$, then we have
$\bar{q}(t)=\bar{u}(t/T)$ is a
non-constant $T$-periodic solution of
{\rm(\ref{eq41})-(\ref{eq42})}.
\end{lemma}

We define a weakly closed subset of $H^1$,
\begin{equation}\label{eq46}
F:=\{u\in E:\int_0^1[V(u)+\frac{1}{2}\min\limits_{y\in \partial V(u)}\langle y, u\rangle] dt=h\}.
\end{equation}

\begin{lemma}\label{lemma3}
If $(V_2)-(V_3)$ hold, then $F\not=\emptyset$.
\end{lemma}

\begin{proof}:
Take $u\in E$ satisfying $\min\limits_{t\in [0,1]}|u(t)|>0$. By condition $(V_2)$, we know $V(0)\leq \mu_2/\mu_1$. We define $$g_u(a):=g(au)=\int_0^1[V(au)+\frac{1}{2}\min_{x\in \partial V(au)}\langle x, au\rangle]dt.$$ Then we have $g_u(0)=g(0)=V(0)\leq\frac{\mu_2}{\mu_1}. \label{eq47}$

We use $(V_2)-(V_3)$ to get
\begin{eqnarray}
g_u(a)&=&g(au)=\int_0^1[V(au)+\frac{1}{2}\min\limits_{x\in \partial V(au)}\langle x, au\rangle]dt\label{eq48}\\
&\geq&(1+\frac{\mu_1}{2})\int_0^1V(au)dt-\frac{\mu_1}{2}\label{eq49}\\
&\rightarrow&+\infty, \ \ as \ \ a\rightarrow +\infty\label{eq50}
\end{eqnarray}
Hence $\forall h>\frac{\mu_2}{\mu_1}$, we know there is $a(u)>0$ such that $a(u)u\in F$.
\end{proof}

\begin{lemma}\label{lemma4}
If $(V_1) - (V_3)$ hold, then for any given $c>0$,
 $f(u)$ satisfies $(CPS)_{F,c;\delta}$
 condition; that is,
 if $\{u_n\}\subset E$ satisfies
\begin{equation}\label{eq51}
dist_\delta(u_n,F)\rightarrow 0,
 f(u_n)\rightarrow c>0,\ \ \ \
(1+\|u_n\|)\min_{y^*\in \partial f(u_n)}\|y^*\|\rightarrow 0,
\end{equation}
then $\{u_n\}$ has a strongly convergent subsequence.
\end{lemma}

\begin{proof}:
 Notice that $\forall u\in E, \int_0^1 u(t)dt=0$; hence,
we know $\|u\|_E:= (\int_0^1|\dot{u}|^2dt)^{1/2}$ is an
equivalent norm on $E$. By $f(u_n)\rightarrow c$, we have
\begin{equation}\label{eq52}
-\frac{1}{2}\|u_n\|^2_E\cdot\int^1_0V(u_n)dt\rightarrow
c-\frac{h}{2}\|u_n\|_E^2.
\end{equation}

By $(V_2)$ we know that $\forall y^*\in \partial f(u_n)$, $\forall x\in \partial V(u_n)$,

\begin{eqnarray}
\langle y^*,u_n\rangle &=&\|u_n\|_E^2\cdot\int^1_0[h-V(u_n)-\frac{1}{2}\langle x,u_n\rangle]dt\nonumber\\
&\leq&\|u_n\|_E^2\int^1_0[h+\frac{\mu_2}{2}-(1+\frac{\mu_1}{2})V(u_n)]dt.\label{eq53}
\end{eqnarray}
By (\ref{eq52}) and (\ref{eq53}) we have
\begin{eqnarray}
\langle y^*,u_n\rangle &\leq&(h+\frac{\mu_2}{2})\|u_n\|_E^2+(1+\frac{\mu_1}{2})(2c-h\|u_n\|_E^2)\nonumber\\
&=&(-\frac{\mu_1}{2}h+\frac{\mu_2}{2})\|u_n\|_E^2+\alpha\label{eq54}
\end{eqnarray}
where $\alpha=2(1+\frac{\mu_1}{2})c$.

Since $h>\frac{\mu_2}{\mu_1}$, then (\ref{eq51}) and (\ref{eq54}) imply $\|u_n\|_{E}$
is bounded.

The rest of the argument to show $\{u_n\}$ has a strongly convergent subsequence
is standard.
\end{proof}
\begin{lemma}\label{lemma5}
Let
\begin{equation}\label{eq55}
G=\{u\in E:\int_0^1[V(u)+\frac{1}{2}\min_{y\in \partial V(u)}\langle y, u\rangle]dt<h\}.
\end{equation}
Then\\
(i). $F$ is the boundary of $G$.\\
(ii). If $(V_1)$ holds, then $F$ is symmetric with respect to
the origin $0$.\\
(iii). If $V(0)<h$ holds, then $0\in G$.
\end{lemma}
It's not difficult to prove the following two Lemmas:
\begin{lemma}\label{lemma6}
 $f(u)$ is weakly lower semi-continuous on $F.$
\end{lemma}
\begin{lemma}\label{lemma7}
$F$ is weakly closed subsets in $H^1$.
\end{lemma}
\begin{lemma}\label{lemma8}
 The functional $f(u)$ has
positive lower bound on $F$.
\end{lemma}
\begin{proof}:
 By the definitions of $f(u)$ and $F$, we have
\begin{equation}\label{eq56}
f(u)=\frac{1}{4}\int^1_0|\dot{u}|^2dt\int^1_0\min_{y\in \partial V(u)}\langle y, u\rangle dt, \ \ u\in
F.
\end{equation}
For $u\in F$ and $(V_2)$, we have
\begin{equation}\label{eq57}
\frac{1}{2}\int^1_0\min_{y\in \partial V(u)}\langle y, u\rangle dt
=\int^1_0 [ h-V(u)]dt\geq \int^1_0 [h-\frac{1}{\mu_1}\min_{y\in \partial V(u)}\langle y, u\rangle-\frac{\mu_2}{\mu_1}]dt,
\end{equation}
\begin{equation}\label{eq58}
\int^1_0\min_{y\in \partial V(u)}\langle y, u\rangle dt \geq\frac{h-\frac{\mu_2}{\mu_1}}{\frac{1}{2}+\frac{1}{\mu_1}}>0.
\end{equation}
So we have the functional $f(u)\geq 0$. Furthermore, we claim that $\inf f(u)>0\label{eq58}$; otherwise, $u(t)=const$ attains the infimum 0.

 If $u\in F$, then by the
symmetry $u(t+1/2)=-u(t)$ or $u(-t)=-u(t)$, we know $u(t)=0,\forall t$. By ($V_2$)
we have $V(0)\leq\frac{\mu_2}{\mu_1}$, by $h>\frac{\mu_2}{\mu_1}$ we get
$V(0)<h$.
From the definition of $F$, $0\notin F$. So $\inf_{F} f(u)>0.\label{eq60}$
Now by Lemmas \ref{lemma6}-\ref{lemma8}, we know $f(u)$ attains the
infimum on $F$, and the minimizer is nonconstant.
\end{proof}
\begin{lemma}\label{lemma9}
$\exists z_1\in H^1$ such that $z_1 \not=0$ and $f(z_1)\leq 0.$
\end{lemma}
\begin{proof}:
 For any given $y_1\not= const$, $\dot{y}_1\not=0$, so
 $\min|\dot{y}_1(t)|>0$.
Let $z_1(t)=Ry_1(t)$, then when $R$ is large enough,
 by condition $(V_3)$ we have
\begin{equation}\label{eq61}
\int_0^1(h-V(z_1))dt\leq0;
\end{equation}
that is, $f(z_1)\leq 0.\label{eq62}$
\end{proof}
\begin{lemma}\label{lemma10}
 $f(0)=0.$
\end{lemma}
\begin{lemma}\label{lemma11}
 $F$ separates $z_1$ and $0$.
\end{lemma}
\begin{proof}:
 By $V(0)<h$, we have that $0\in G$.
By $(V_2)$ and $(V_3)$ and $h>\frac{\mu_2}{\mu_1}$, we can choose R large enough
such that
\begin{eqnarray}
 z_1&=&Ry_1\in \{u\in H^1:\int_0^1[V(u)+\frac{1}{2}\min_{y\in \partial V(u)}\langle y, u\rangle]dt\label{eq63}\\ &\geq&(1+\frac{\mu_1}{2})\int_0^1V(u)dt-\frac{\mu_1}{2}\label{eq64}\\
 &\geq&(1+\frac{\mu_1}{2})h-\frac{\mu_1}{2}>h\}.\label{eq65}
\end{eqnarray}
So $F$ separates $z_1$ and $0$.
\end{proof}

Theorem \ref{thm 9} now follows from Theorem \ref{thm 7}

\section{Conclusions}

Since Ekeland's variational principle imposes less restriction on the functional, we found it very useful in proving our Generalized Mountain Pass Lemma with weaker assumptions. We were able to establish an immediate application for our Generalized Mountain Pass Lemma to Hamiltonian systems  with Lipschtitz potential and a fixed energy. It would be interesting to see what role it can play for other differential equations.

\section*{Acknowledgements}

This research was partially supported by NSF of China(No.11671278) and the Grant for the Advisors of Ph.D students(No.20120181110060).



\end{document}